# Bounds of Shannon entropy and Extropy and their application in exploring the extreme value behavior of a large set of data

**Konstantinos Zografos**[0]

University of Ioannina, Department of Mathematics, 451 10 Ioannina, Greece

## Abstract

This paper derives bounds for two omnipresent information theoretic measures, the Shannon entropy and its complementary dual, the extropy. Based on a large size data set from a log-concave model, the said bounds are obtained for the entropy and the extropy of the distribution of the largest order statistic and the respective normalized sequence, in the extreme value theory setting. A characterization of the exponential distribution is provided as the model that maximizes the Shannon entropy and the extropy which are associated with the distribution of the maximum value, in a large sample size regime. This characterization is exploited to provide an alternative, immediate proof of the convergence of Shannon entropy and extropy of the normalized maxima of a large size sample to the respective measures for the Gumbel distribution, studied recently for Shannon entropy in Johnson (2024) and references therein.

## 1    Introduction

Having in hand a random sample from a probabilistic model, the largest/smallest order statistics play a fundamental role in statistical inference, especially when the domain of the probabilistic model that drives the data depends on an unknown parameter of interest. On the other hand, order statistics in general are the fundamental tools in reliability as they are usually related with the failure time of a tested item. Moreover, the largest order statistic plays a fundamental role in extreme value theory and therefore it is the basis in developing the extreme value distributions in cases where the size of the random sample approaches to infinity. On the other hand, Shannon entropy of a continuous distribution is an omnipresent, global quantity. Shannon entropy expresses the degree of uncertainty about the outcome of a random experiment. Large values of Shannon entropy are associated with high uncertainty or low information gain which can be thought of as the purity in a system (cf. also Zaid et al. (2022), p. 13). In this context, it would by maybe of interest to study the behavior of Shannon entropy when it is applied to the distribution of the largest order statistic in the environment of the extreme value theory, that is, in the case of an infinite number of independent and identically distributed observations. In this frame, this paper is focused in the investigation of Shannon entropy and its dual respective, namely the extropy, when these measures are applied to the distribution of the largest order statistic.

To formulate the work, consider a random sample $X_1, X_2, ..., X_n$ of size $n$ from a continuous probability distribution function $F$ with respective probability density function $f$. The Shannon entropy associated to the distribution with density $f$ is defined by

$$\mathcal{H}(X) = -\int_{\mathbb{R}} f(x) \ln f(x) dx. \tag{1}$$

---
[0]Corresponding author e-mail address: kzograf@uoi.gr
ORCID iD: https://orcid.org/0000-0003-3677-7850

It is a well known and broadly applied quantity with applications in almost each field of science and engineering. On the other hand, the extropy of a random variable or the extropy of the respective density function has received considerable attention the last decade in the field of statistical information theory. Based on Zografos (2023), p. 296 and references appeared therein, the extropy of a random variable $X$ with an absolutely continuous distribution function $F$ and respective probability density function $f$ is defined by

$$\mathcal{J}(X) = -\frac{1}{2} \int_{\mathbb{R}} f^2(x)dx, \tag{2}$$

and it has been introduced in the statistical literature as the complement of Shannon entropy. This measure is directly connected with other informational quantities, like Onicescu information energy or Golomb information function, discussed in Zografos (2023). More information about extropy can be extracted from the papers by Lad et al. (2015), Qiu (2017), Toomaj et al. (2023) and the references appeared therein.

Based on the random sample $X_1, X_2, ..., X_n$, from a continuous probability distribution function $F$ with respective probability density function $f$, let $X_{(1)} \leq X_{(2)} \leq ... \leq X_{(n)}$ be the ordered sample with $X_{(n)} = \max\{X_1, X_2, ..., X_n\}$ to be the largest order statistics. $X_{(n)}$ is a fundamental quantity in statistical inference and in reliability theory because it quantifies the largest failure time of a set of items being under a test. It is well known that the density of $X_{(n)}$ is given by

$$f_{X_{(n)}}(x) = nF^{n-1}(x)f(x), \tag{3}$$

and it can be used in order to obtain the entropy and the extropy associated to the distribution of $X_{(n)}$, say $\mathcal{H}(X_{(n)})$ and $\mathcal{J}(X_{(n)})$. Due to the importance of $X_{(n)}$ in statistical practice, several papers have been devoted to the study of Shannon entropy and extropy of order statistics. A list of references of the existing literature in this subject is given in the paper by Qiu (2017), p. 52, where some characterization results and bounds of extropy of order statistics and record values are presented (cf. Qiu and Jia (2018)). Several properties of the information generating function of $k$-record values are studied in Chacko and Grace (2023) while the recent paper by Mohammadi and Hashmpour (2023) presents and studies an extropy based inaccuracy measure of order statistiscs. All the above mentioned papers and the references appeared therein are gathered the existing literature on some measures of information and order statistics and records.

Although the distribution of the maximum of $n$, $n < \infty$, independent and identically distributed (i.i.d.) random variables is well defined in (3), the limit of $X_{(n)}$ may degenerate when $n \to \infty$. To avoid this degeneracy, we consider the sequence $\frac{X_{(n)} - b_n}{a_n}$, for some suitable sequences $a_n$ and $b_n, n \geq 1$, with $a_n > 0$. In this direction, it is well known from the area of the extreme value distributions that $\frac{X_{(n)} - b_n}{a_n}$ converges in distribution to some non-degenerate limiting distribution $G$ if and only if the distribution $F$ is in the domain of attraction of the distribution $G$, something which is denoted by $F \in \mathcal{D}(G)$ (cf. Beirlant et al. 2004, p. 46-47, or Kotz and Nadarajah, 2000, p. 7-9). The study of the behavior of Shannon entropy, extropy and other information measures of the largest order statistic $X_{(n)}$, when $n \to \infty$, has been considered in the recent literature. Recent advances on this topic constitute the papers by Zaid et al. (2022) where the interest is focused on some measures of information when they are applying in the $q$-gereralized extreme value distributions. The work by Saeb (2023), cf. also Ravi and Saeb (2014), concentrates on the convergence of Shannon entropy of the density function of the normalized maxima of a random sample to the Shannon entropy of the limit max-stable low and it moreover shows that this convergence is equivalent to the convergence of the relative entropy to zero, under some conditions.

The recent paper by Johnson (2024) introduces a type of score function, namely the max-score function, which is used to prove convergence to the Gumbel distribution in the strong sense of relative entropy. Based on the introduction of this last paper, it provides, in contrast with the work by Ravi and Saeb (2014) and Saeb (2023), "an elementary and direct proof under relatively simple conditions which hopefully gives some insight into why convergence to the Gumbel takes place, rather than to necessarily provide the strongest possible result.".

Based on the above exposition, it seems natural to be interested in the study of the behavior of $\mathcal{H}(X_{(n)})$ and $\mathcal{J}(X_{(n)})$, and their respective normalized versions, when $n \to \infty$, as it is natural in extreme value theory to be interested in finding the possible limit distributions of $X_{(n)}$ (cf. Beirlant et al. 2004, p. 46). In this context, bounds for the Shannon entropy and the extropy of $X_{(n)}$ will be derived in the next section when the independent identically distributed (i.i.d.) random variables $X_1, X_2, ..., X_n$ are governed by a log-concave distribution. A characterization of the exponential distribution will be provided and it will be shown that it is the model which maximizes $\mathcal{H}(X_{(n)})$ and $\mathcal{J}(X_{(n)})$, when $n \to \infty$. The limiting behavior of $\mathcal{H}$ and $\mathcal{J}$ of the normalized maxima will be studied in the subsequent Section 3 and a similar characterization of the exponential distribution will be also provided. Then, the convergence of $\mathcal{H}(X_{(n)})$ and $\mathcal{J}(X_{(n)})$ to the respective Shannon entropy and the extropy of the Gumbel distribution, discussed in Ravi and Saeb (2014), Saeb (2023) and Johnson (2024), will be implicitly derived, from another point of view. The paper is completed with some conclusions and an Appendix, the Section 5, of the proofs of the statements. Hence, in summary, this paper studies the behavior of the Shannon entropy and the extropy in the extreme value regime, providing with bounds for these measures, when $n \to \infty$, inside the broad family of log-concave distributions. A subsequent characterization of the exponential distribution is also provided.

## 2 Bounds for Shannon entropy and Extropy

This section is devoted to the study of two omnipresent in science and engineering measures, namely, the Shannon entropy and the extropy of a probability distribution. More precisely the interest is focused in the investigation of the behavior of these indices when they are applying to the largest order statistic of a set of independent random variables which are governed by the same probability distribution which is assumed to be log-concave.

### 2.1 Shannon entropy of the largest order statistic

Let a random sample $X_1, X_2, ..., X_n$ of size $n$ from a continuous probability distribution function $F$ with respective probability density function $f$. Let also $X_{(1)} \leq X_{(2)} \leq ... \leq X_{(n)}$ be the respective order statistics with $X_{(n)} = \max\{X_1, X_2, ..., X_n\}$. The Shannon entropy associated to the distribution with density $f$ is defined by (1). An application of (1) in the specific case (3) of the distribution of $X_{(n)}$ leads to the Shannon entropy of $X_{(n)}$ which can be also obtained by applying Corollary 1 of Zografos and Balakrishnan (2009), for $\alpha = n$ and $\beta = 1$, and then the Shannon entropy of $X_{(n)} = \max\{X_1, X_2, ..., X_n\}$ is given by,

$$\mathcal{H}(X_{(n)}) = \ln B(n, 1) - (n-1)[\Psi(n) - \Psi(n+1)] - E_Y[\ln f(F^{-1}(Y))],$$

where $Y$ is a random variable which follows a beta distribution $Y \sim Beta(n, 1)$ and $B(\alpha, \beta) = \Gamma(\alpha)\Gamma(\beta)/\Gamma(\alpha + \beta)$, with $B$ and $\Gamma$ being the beta and gamma functions, respectively. In this

context, if we denote by

$$I(Y) = f(F^{-1}(Y)), \ Y \sim Beta(n,1), \tag{4}$$

then $\mathcal{H}(X_{(n)})$ is simplified as follows,

$$\mathcal{H}(X_{(n)}) = 1 - \ln n - \frac{1}{n} - E_Y[\ln I(Y)), \text{ or} \tag{5}$$

$$\mathcal{H}(X_{(n)}) = 1 - \ln n - \frac{1}{n} - \int\limits_{0}^{1} ny^{n-1} \ln f(F^{-1}(y))dy, \tag{6}$$

taking into account that for the digamma function $\Psi$ it holds that $\Psi(n+1) - \Psi(n) = \frac{1}{n}$. $\mathcal{H}(X_{(n)})$ in (5) or (6) formulates the Shannon entropy of the extreme value distribution of the maximum $X_{(n)}$ of a sequence of i.i.d. data $X_1, X_2, ..., X_n$.

A fundamental role in the derivations which will be follow, plays the real function

$$I(t) = f(F^{-1}(t)), \ 0 < t < 1, \tag{7}$$

which is motivated from (4) and it is named the *density quantile function* by David and Nagaraja (2003), p. 84. It has been investigated and exploited in a series of publications by Bobkov and his colleagues, cf. for example, Bobkov (1996, 1999), Bobkov and Madiman (2011), Bobkov and Ledoux (2019). If the density $f$, included in its definition, is log-concave, then the function $I$, which is defined by (7), obeys nice properties, like positiveness and concavity. Based on the said properties of $I$, next theorem formulates some bounds for the Shannon entropy $\mathcal{H}(X_{(n)})$. In the statement of the theorem, some well known results from Abramowitz and Stegun (1970), p. 255, 258, Mathai (1993), p. 4, 11, among others, will be used. In this frame,

$$\gamma = \lim_{n \to \infty} (H_n - \ln n) = \lim_{n \to \infty} \left( \sum_{k=1}^{n} \frac{1}{k} - \ln n \right), \tag{8}$$

is the Euler's constant while $H_n$ is the $n$-th harmonic number, defined by

$$H_n = \sum_{k=1}^{n} \frac{1}{k}. \tag{9}$$

**Theorem 1** *Let a random sample $X_1, X_2, ..., X_n$ of size $n$ from a continuous probability distribution function $F$ with respective probability density function $f$, supported on some interval $(\alpha, \beta)$, finite or not, where $f$ is positive and $\ln f$ is concave. Then,*
*(i) It holds that,*

$$1 - \ln n - \frac{1}{n} \leq \mathcal{H}(X_{(n)}) \leq 1 - \ln \left[ 2I \left( \frac{1}{2} \right) \right] - \ln n - \frac{1}{n} + \ln 2 + \sum_{k=1}^{n} \frac{1}{k} - \sum_{k=1}^{n-1} \frac{1}{k2^k},$$

*with $I$ defined by (7).*
*(ii) The limiting value of $\mathcal{H}(X_{(n)})$, as $n \to \infty$, obeys the next inequality*

$$-\infty \leq \lim_{n \to \infty} \mathcal{H}(X_{(n)}) \leq 1 - \ln \left[ 2I \left( \frac{1}{2} \right) \right] + \gamma,$$

*where $\gamma$ is the Euler constant, defined in (8) and (9).*

The proof of the Theorem is presented in subsection 5.1 of the proofs of the theoretical results.

It is obtained in Table 1 the exact form of $\mathcal{H}(X_{(n)})$ and its limiting value $\lim_{n\to\infty}\mathcal{H}(X_{(n)})$ for specific distributions. The explicit form of the function $I(t)$, defined by (7), is also given for each one of these distributions. The last column includes in addition the upper bound of the $\lim_{n\to\infty}\mathcal{H}(X_{(n)})$, obtained in part (ii) of the previous theorem and it is denoted by

$$\mathcal{U}\mathcal{B}^{\mathcal{H}} = 1 - \ln\left[2I\left(\frac{1}{2}\right)\right] + \gamma. \tag{10}$$

A subscript is used in the notation of $\mathcal{U}\mathcal{B}^{\mathcal{H}}$, in (10), with the first letter of the distribution considered, each time. The derivations of the formulas is the result of standard algebraic manipulations which take also into account results from Zografos and Balakrishnan (2009), p. 348-350, and standard properties of special functions like that of

$$\Psi(n+1) - \Psi(n) = \frac{1}{n} = [H_n - \gamma] - [H_{n-1} - \gamma],$$

where $\Psi$ is the digamma function, $H_n$ is the $n$-th harmonic number, defined by (9) and $\gamma$ is the Euler's constant, related with $H_n$ by $\gamma = \lim_{n\to\infty}(H_n - \ln n) = 0.57721...$, given in (8), (cf. Mathai (1993), p. 11, (1.4.6), among others).

| | **Distribution and $I(t)$** | $\mathcal{H}(X_{(n)})$, $\lim_{n\to\infty}\mathcal{H}(X_{(n)})$ **and upper bound** $\mathcal{U}\mathcal{B}^{\mathcal{H}}$ **in (10)** |
|---|---|---|
| **1** | Uniform, $U(0,\theta)$ <br> $f(x) = 1/\theta,$ <br> $F(x) = x/\theta, 0 < x < \theta$ <br> $I(t) = 1/\theta, 0 < t < 1$ | $\mathcal{H}(X_{(n)}) = 1 - \ln n - \frac{1}{n} + \ln\theta$ <br> $\lim_{n\to\infty}\mathcal{H}(X_{(n)}) = -\infty$ <br> $\mathcal{U}\mathcal{B}_U^{\mathcal{H}} = 1 - \ln 2 + \ln\theta + \gamma$ |
| **2** | Exponential, $Exp(\theta)$ <br> $f(x) = \theta e^{-\theta x},$ <br> $F(x) = 1 - e^{-\theta x}, x > 0, \theta > 0$ <br> $I(t) = \theta(1-t), 0 < t < 1$ | $\mathcal{H}(X_{(n)}) = 1 - \ln n - \frac{1}{n} - \ln\theta + \sum_{k=1}^{n}\frac{1}{k}$ <br> $\lim_{n\to\infty}\mathcal{H}(X_{(n)}) = 1 - \ln\theta + \gamma, \ \gamma = 0.57721...$ <br> $\mathcal{U}\mathcal{B}_{\exp}^{\mathcal{H}} = 1 - \ln\theta + \gamma$ |
| **3** | Logistic <br> $f(x) = \frac{\theta e^{-\theta x}}{(1+e^{-\theta x})^2},$ <br> $F(x) = \frac{1}{1+e^{-\theta x}}, x \in \mathbb{R}, \theta > 0$ <br> $I(t) = \theta t(1-t), 0 < t < 1$ | $\mathcal{H}(X_{(n)}) = 1 - \ln n - \ln\theta + \sum_{k=1}^{n}\frac{1}{k}$ <br> $\lim_{n\to\infty}\mathcal{H}(X_{(n)}) = 1 - \ln\theta + \gamma, \ \gamma = 0.57721...$ <br> $\mathcal{U}\mathcal{B}_{\log}^{\mathcal{H}} = 1 + \ln 2 - \ln\theta + \gamma$ |
| **4** | Pareto <br> $f(x) = \frac{\nu\theta^{\nu}}{x^{\nu+1}},$ <br> $F(x) = 1 - \left(\frac{\theta}{x}\right)^{\nu}, x \geq \theta > 0, \nu > 0$ <br> $I(t) = \frac{\nu(1-t)^{(\nu+1)/\nu}}{\theta}, 0 < t < 1$ | $\mathcal{H}(X_{(n)}) = 1 + \frac{1}{\nu}\ln n - \frac{1}{n} - \ln\frac{\nu}{\theta} + \frac{\nu+1}{\nu}\left(\sum_{k=1}^{n}\frac{1}{k} - \ln n\right)$ <br> $\lim_{n\to\infty}\mathcal{H}(X_{(n)}) = +\infty$ <br> $\mathcal{U}\mathcal{B}_P^{\mathcal{H}} = 1 + \frac{1}{\nu}\ln 2 - \ln\nu + \ln\theta + \gamma$ |
| **5** | Power-function <br> $f(x) = \nu\theta^{\nu}x^{\nu-1}$ <br> $F(x) = (\theta x)^{\nu}, 0 < x < \frac{1}{\theta}, \nu > 0$ <br> $I(t) = \nu\theta t^{(\nu-1)/\nu}, 0 < t < 1$ | $\mathcal{H}(X_{(n)}) = 1 - \ln n - \ln(\nu\theta) - \frac{1}{\nu n}$ <br> $\lim_{n\to\infty}\mathcal{H}(X_{(n)}) = -\infty$ <br> $\mathcal{U}\mathcal{B}_{PF}^{\mathcal{H}} = 1 - \frac{1}{\nu}\ln 2 - \ln\nu - \ln\theta + \gamma$ |

**Table 1:** Analytic Expressions for $\mathcal{H}(X_{(n)})$, $\lim_{n\to\infty}\mathcal{H}(X_{(n)})$ and $\mathcal{U}\mathcal{B}^{\mathcal{H}}$ for specific $F$ with respective density $f$

It should be mentioned, at this point, that the expression for $\mathcal{H}(X_{(n)})$ of the exponential distribution can be also obtained by using a similar expression, derived on p. 693 of Nadarajah and Kotz (2006), for the beta exponential distribution.

**Remark 1** *We observe that in all the cases the inequality* $\lim_{n\to\infty} \mathcal{H}(X_{(n)}) \leq \mathcal{UB}^{\mathcal{H}}$*, which is stated in part (ii) of Theorem 1, is valid, except of the case of Pareto distribution. This is well expected as the Pareto distribution is not log-concave and hence the basic assumption of log-concavity of Theorem 1 is violated for this distribution. Moreover, we observe that* $\lim_{n\to\infty} \mathcal{H}(X_{(n)})$ *attains its maximum value* $\mathcal{UB}^{\mathcal{H}}_{\exp}$ *only in the case of the exponential distribution. This point motivates the investigation of a characterization of the exponential distribution in terms of the attainance of the upper bound (10) by the limit of Shannon entropy of the largest order statistic. A similar characterization of an exponential type distribution has been derived in Chacko and Grace (2023) by means of the information generating function of records.*

The characterization of the exponential distribution, mentioned in the previous remark, is formulated in the next theorem.

**Theorem 2** *Let a random sample* $X_1, X_2, ..., X_n$ *of size* $n$ *from a continuous probability distribution function* $F_\theta$ *with respective probability density function* $f_\theta$*, supported on some interval* $(\alpha, \beta)$*, finite or not, where* $f_\theta$ *is positive and* $\ln f_\theta$ *is concave. Let moreover that* $f_\theta$ *depends on an real parameter* $\theta$*. Then, the* $\lim_{n\to\infty} \mathcal{H}(X_{(n)})$ *attains its maximum value* $\mathcal{UB}^{\mathcal{H}}$*, given in (10), if and only if* $F_\theta$ *is the distribution function of the exponential distribution.*

The proof of the theorem is given in subsection 5.2 of the proofs of the theoretic results. Intuitively speaking, the above characterization concludes that the exponential distribution maximizes the uncertainty, as it is quantified by Shannon entropy, which is associated with the maximum value in a sample of size $n$, when $n \to \infty$.

Figure 1 illustrates Theorems 1 and 2 when the data are coming from the standard exponential distribution with density $e^{-x}, x > 0$.

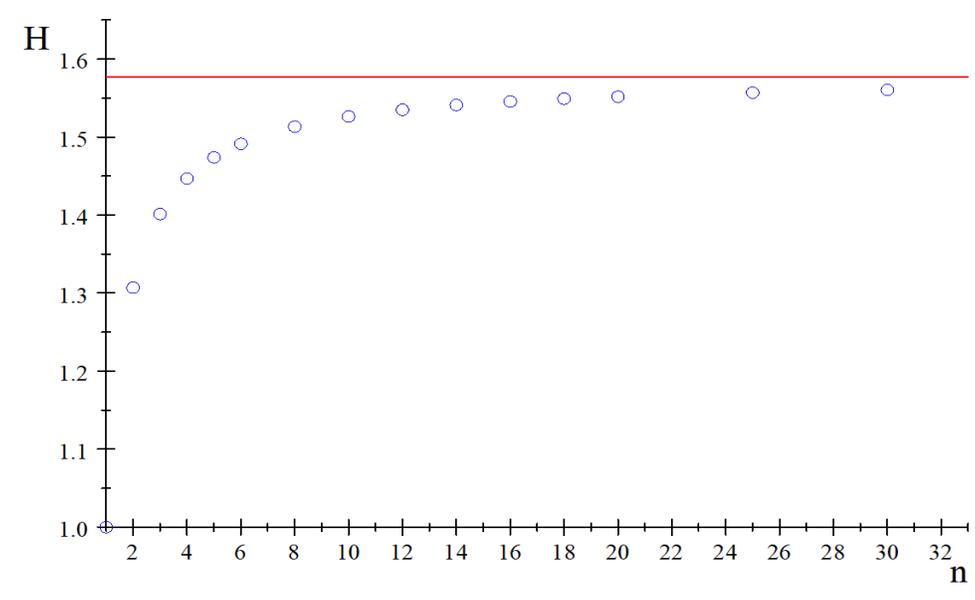

Figure 1: Plot of $\mathcal{H}(X_{(n)})$ (blue points/cycles) and $\mathcal{UB}^{\mathcal{H}}_{\exp} = 1 + \gamma$ (red solid) for different values of $n$ and the standard exponential distribution $e^{-x}, x > 0$.

## 2.2  Extropy of the largest order statistic

The extropy of a random variable or the extropy of the respective density function has been defined by (2). It has been studied recently by Lad et al. (2015), and based on this paper, the extropy of a random variable $X$ with an absolutely continuous distribution function $F$ and respective probability density function $f$ is given by

$$\mathcal{J}(X) = -\frac{1}{2}\int_{\mathbb{R}} f^2(x)d\mu.$$

It has been introduced in the statistical literature as the complementary dual of Shannon entropy, underlined in the title of the paper by Lad et al. (2015).

Let again a random sample $X_1, X_2, ..., X_n$ of size $n$ from a continuous probability distribution function $F$ with respective probability density function $f$ and let $X_{(n)} = \max\{X_1, X_2, ..., X_n\}$ be the largest order statistics. It is well known from (3) that the density of $X_{(n)}$ is given by

$$f_{X_{(n)}}(x) = nF^{n-1}(x)f(x),$$

and then its extropy,

$$\mathcal{J}(X_{(n)}) = -\frac{1}{2}\int_{\mathbb{R}} f^2_{X_{(n)}}(x)dx = -\frac{n^2}{2}\int_{\mathbb{R}} F^{2n-2}(x)f^2(x)dx. \tag{11}$$

Assuming that $F$ is strictly increasing and using the transformation $t = F(x)$, $0 < t < 1$, $x = F^{-1}(t)$, where $F^{-1}$ is the inverse of $F$ and $\frac{dx}{dt} = \frac{d}{dt}F^{-1}(t) = \frac{1}{F'(F^{-1}(t))} = \frac{1}{f(F^{-1}(t))}$, (11) leads to

$$\mathcal{J}(X_{(n)}) = -\frac{n^2}{2}\int_0^1 t^{2n-2}\left[f\left(F^{-1}(t)\right)\right]^2\frac{1}{f\left(F^{-1}(t)\right)}dt = -\frac{n^2}{2}\int_0^1 t^{2n-2}f\left(F^{-1}(t)\right)dt,$$

or

$$\mathcal{J}(X_{(n)}) = -\frac{n^2}{2}\int_0^1 t^{2n-2}I(t)dt, \tag{12}$$

where $I(t) = f\left(F^{-1}(t)\right), 0 < t < 1$, is defined by (7). Hence, the following lemma has been obtained and it formulates the extropy of the largest order statistic.

**Lemma 1** *Let a random sample $X_1, X_2, ..., X_n$ of size $n$ from a strictly increasing continuous probability distribution function $F$ with respective probability density function $f$. Let also $X_{(1)} \leq X_{(2)} \leq ... \leq X_{(n)}$ be the respective order statistics with $X_{(n)} = \max\{X_1, X_2, ..., X_n\}$. Then the extropy of $X_{(n)}$ is given by*

$$\mathcal{J}(X_{(n)}) = -\frac{n^2}{2}\int_0^1 t^{2n-2}I(t)dt,$$

*where $I(t) = f\left(F^{-1}(t)\right), 0 < t < 1$, is defined by (7).*

A question is raised at this point, is it possible to derive bounds for the extropy $\mathcal{J}(X_{(n)})$, similar to those of Theorem 1? An answer is provided in the next theorem, the proof of which is sketched in subsection 5.3 of the proofs of the statements. Lower bounds for the extropy of order statistics and record values have been obtained in Qiu (2017) in a different setting.

**Theorem 3** *Let a random sample $X_1, X_2, ..., X_n$ of size $n$ from a continuous probability distribution function $F$ with respective probability density function $f$, supported on some interval $(\alpha, \beta)$, finite or not, where $f$ is positive and $\ln f$ is concave. Then,*
*(i) it holds that,*

$$-\frac{n}{2}I\left(\frac{1}{2}\right) \leq \mathcal{J}(X_{(n)}) \leq -n^2 I\left(\frac{1}{2}\right)\left(\frac{1}{2n-1} - \frac{1}{2n} - \frac{1}{(2n-1)2^{2n-1}} + \frac{2}{2n2^{2n}}\right),$$

*with $I$ defined by (7).*
*(ii) The limiting value of $\mathcal{J}(X_{(n)})$, as $n \to \infty$, obeys the next inequality*

$$-\infty \leq \lim_{n\to\infty} \mathcal{J}(X_{(n)}) \leq -\frac{1}{4}I\left(\frac{1}{2}\right).$$

The analog of Table 1 for the extropy is the following. It is given in Table 2 the exact form of $\mathcal{J}(X_{(n)})$ and its limiting value $\lim_{n\to\infty}\mathcal{J}(X_{(n)})$ for specific distributions. The explicit form of the function $I(t)$, defined by (7), is also given for each one of these distributions while the last column includes in addition the upper bound of the $\lim_{n\to\infty}\mathcal{J}(X_{(n)})$, obtained in part (ii) of the previous theorem and it is denoted by

$$\mathcal{UB}^{\mathcal{J}} = -\frac{1}{4}I\left(\frac{1}{2}\right). \tag{13}$$

| | **Distribution and $I(t)$** | $\mathcal{J}(X_{(n)})$, $\lim_{n\to\infty}\mathcal{J}(X_{(n)})$ **and upper bound** $\mathcal{UB}^{\mathcal{J}}$ **in (13)** |
|---|---|---|
| **1** | Uniform, $U(0,\theta)$ <br> $f(x) = 1/\theta$, <br> $F(x) = x/\theta, 0 < x < \theta$ <br> $I(t) = 1/\theta, 0 < t < 1$ | $\mathcal{J}(X_{(n)}) = -\frac{n^2}{2(2n-1)\theta}$ <br> $\lim_{n\to\infty}\mathcal{J}(X_{(n)}) = -\infty$ <br> $\mathcal{UB}^{\mathcal{J}}_U = -\frac{1}{4\theta}$ |
| **2** | Exponential, $Exp(\theta)$ <br> $f(x) = \theta e^{-\theta x}$, <br> $F(x) = 1 - e^{-\theta x}, x > 0, \theta > 0$ <br> $I(t) = \theta(1-t), 0 < t < 1$ | $\mathcal{J}(X_{(n)}) = -\frac{n^2\theta}{2(2n-1)} + \frac{n^2\theta}{4n}$ <br> $\lim_{n\to\infty}\mathcal{J}(X_{(n)}) = -\frac{\theta}{8}$ <br> $\mathcal{UB}^{\mathcal{J}}_{\exp} = -\frac{\theta}{8}$ |
| **3** | Logistic <br> $f(x) = \frac{\theta e^{-\theta x}}{(1+e^{-\theta x})^2}$, <br> $F(x) = \frac{1}{1+e^{-\theta x}}, x \in \mathbb{R}, \theta > 0$ <br> $I(t) = \theta t(1-t), 0 < t < 1$ | $\mathcal{J}(X_{(n)}) = -\frac{n\theta}{4} + \frac{n^2\theta}{2(2n+1)}$ <br> $\lim_{n\to\infty}\mathcal{J}(X_{(n)}) = -\frac{\theta}{8}$ <br> $\mathcal{UB}^{\mathcal{J}}_{\log} = -\frac{\theta}{16}$ |
| **4** | Pareto <br> $f(x) = \frac{\nu\theta^\nu}{x^{\nu+1}}$, <br> $F(x) = 1 - \left(\frac{\theta}{x}\right)^\nu, x \geq \theta > 0, \nu > 0$ <br> $I(t) = \frac{\nu(1-t)^{(\nu+1)/\nu}}{\theta}, 0 < t < 1$ | $\mathcal{J}(X_{(n)}) = -\frac{\nu n^2}{2\theta}B\left(2n-1, \frac{2\nu+1}{\nu}\right)$ <br> $\lim_{n\to\infty}\mathcal{J}(X_{(n)}) = 0 \times (-\infty)$ <br> $\mathcal{UB}^{\mathcal{J}}_P = -\frac{\nu}{\theta}\frac{1}{2^{(3\nu+1)/\nu}}$ |
| **5** | Power-function <br> $f(x) = \nu\theta^\nu x^{\nu-1}$ <br> $F(x) = (\theta x)^\nu, 0 < x < \frac{1}{\theta}, \nu > 0$ <br> $I(t) = \nu\theta t^{(\nu-1)/\nu}, 0 < t < 1$ | $\mathcal{J}(X_{(n)}) = -\frac{\nu^2\nu^2\theta}{2(2n\nu-1)}$ <br> $\lim_{n\to\infty}\mathcal{J}(X_{(n)}) = -\infty$ <br> $\mathcal{UB}^{\mathcal{J}}_{PF} = -\frac{\nu\theta}{2^{(3\nu-1)/\nu}}$ |

**Table 2:** Analytic Expressions for $\mathcal{J}(X_{(n)})$, $\lim\limits_{n\to\infty}\mathcal{J}(X_{(n)})$ and $\mathcal{UB}^{\mathcal{J}}$ for specific $F$ with respective density $f$

**Remark 2** *In complete analogy with Remark 1, we observe that in all the cases the inequality which is stated in part (ii) of Theorem 3 is valid, except of the case of Pareto distribution for the same reasons as those which are mentioned in Remark 1. Moreover, we observe that $\lim_{n\to\infty}\mathcal{J}(X_{(n)})$* *attains its maximum value $\mathcal{UB}_{\exp}^{\mathcal{J}}$ only for the exponential distribution. This remark motivates an additional characterization of the exponential distribution in terms of its extropy and the Upper Bound of part (ii) of the previous Theorem 3.*

The next theorem formulates the characterization of the exponential distribution as the model that maximizes the $\lim_{n\to\infty} \mathcal{J}(X_{(n)})$. The proof of this characterization is outlined in subsection 5.4 of the proofs of the statements.

**Theorem 4** *Let a random sample $X_1, X_2, ..., X_n$ of size $n$ from a continuous probability distribution function $F_\theta$ with respective probability density function $f_\theta$, supported on some interval $(\alpha, \beta)$, finite or not, where $f_\theta$ is positive and $\ln f_\theta$ is concave. Let moreover that $f_\theta$ depends on an real parameter $\theta$. Then, the $\lim_{n\to\infty} \mathcal{J}(X_{(n)})$ attains its maximum value $\mathcal{UB}^{\mathcal{J}}$, given in (13), if and only if $F_\theta$ is the distribution function of the exponential distribution.*

## 3 Extreme Value Distribution

We observe in the above tables that there is a diverse behavior of $\lim_{n\to\infty} \mathcal{H}(X_{(n)})$ or $\lim_{n\to\infty} \mathcal{J}(X_{(n)})$. They are finite for some distributions while for some others the same limit does not exist, it is $-\infty$ or $+\infty$. This limiting behavior of $\mathcal{H}(X_{(n)})$ and $\mathcal{J}(X_{(n)})$, when $n \to \infty$, is expected if we take into account that the limit of $X_{(n)}$ may degenerate when $n \to \infty$. To avoid this degeneracy, motivated by the extreme value theory, we consider the sequence $\frac{X_{(n)} - b_n}{a_n}$, for some suitable sequences $a_n$ and $b_n, n \geq 1$, with $a_n > 0$. In this direction, it is well known from the area of the extreme value distributions that $\frac{X_{(n)} - b_n}{a_n}$ converges in distribution to some non-degenerate limiting distribution $G$ if and only if the distribution $F$ is in the domain of attraction of the distribution $G$, something which is denoted by $F \in \mathcal{D}(G)$ (cf. Beirlant et al. 2004, p. 46-47, or Kotz and Nadarajah, 2000, p. 7-9). Hence, it is reasonable to study the limiting behavior of Shannon entropy and extropy of the normalizing sequence $\frac{X_{(n)} - b_n}{a_n}$ by using a parallel exposition like that of the previous section. However, firstly, we will present some prerequisites from the theory of extreme value distributions for the sake of completeness of the presentation.

Let a random sample $X_1, X_2, ..., X_n$ of size $n$ from the continuous probability distribution function $F$. Then, $\frac{X_{(n)} - b_n}{a_n}$ converges in distribution to some non-degenerate limiting distribution $G$ for some sequences $a_n$ and $b_n, n \geq 1$, with $a_n > 0$, if and only if $F \in \mathcal{D}(G)$, that is, the distribution $F$ is in the domain of attraction of the distribution $G$. Therefore, if $F \in \mathcal{D}(G)$ then,

$$\lim_{n\to\infty} P\left(\frac{X_{(n)} - b_n}{a_n} \leq x\right) = G_\xi(x), \tag{14}$$

and based on Beirlant et al. (2004), p. 47, or de Haan and Ferreira (2006), p. 6, if $F \in \mathcal{D}(G)$ then, all extreme value distributions

$$G_\xi(x) = \exp\left(-(1 + \xi x)^{-1/\xi}\right), \ 1 + \xi x > 0, \ \xi \in \mathbb{R}, \tag{15}$$

can occur as limits of (14). For $\xi = 0$, the right-hand side of (15) is interpreted as $\exp\left(-e^{-x}\right)$ while $\xi$ is the extreme value index. Hence, for $\xi = 0$,

$$G_0(x) = \exp\left(-e^{-x}\right), x \in \mathbb{R}. \tag{16}$$

The distribution function in (15) with respective density

$$g_\xi(x) = (1 + \xi x)^{-(\xi+1)/\xi} \exp\left(-(1 + \xi x)^{-1/\xi}\right), \ 1 + \xi x > 0, \ \xi \in \mathbb{R}, \tag{17}$$

is the well known class of extreme value distributions, introduced, according to Kotz and Nadarajah (2000), p. 61, by Jenkinson (1955) and it also referred to as the von Mises type extreme value distributions or the von Mises-Jenkinson type distributions, in view of the previous mentioned book. The parameter $\xi$ is the shape parameter of the model and it may be used to model a wide range of tail behavior (cf. Kotz and Nadarajah, 2000, p. 62 for more details). For $\xi = 0$, the right-hand side of (17) is interpreted as

$$g_0(x) = \exp\left(-(x + e^{-x})\right), x \in \mathbb{R}. \tag{18}$$

The next proposition investigates the log-concavity properties of $g$ in (17) and it moreover provides some results related to the Shannon entropy and the extropy of the largest order statistic from the distribution in (15) or (17). The log-concavity of the extreme value distribution has been studied in Müller and Rufibach (2008). For the sake of completeness of the presentation an elementary proof of the log-concavity and the whole proof of the next proposition is outlined in subsection 5.5 of the proofs of the statements.

**Proposition 1** *Let a random sample* $X_1, X_2, ..., X_n$ *of size n from the continuous probability distribution function* $G_\xi$ *in (15) or with respective probability density function* $g_\xi$ *in (17). Then,*
*(i) The generalized extreme value density function* $g_\xi$*, in (17), is log-concave for* $\xi = 0$ *or* $-1 < \xi < 0$.
*(ii) The function* $I(t) = g_\xi(G_\xi^{-1}(t))$, $0 < t < 1$, *in (7), for this distribution is given by,*

$$I(t) = t(-\ln t)^{\xi+1}, \ 0 < t < 1, \ \xi \in \mathbb{R},$$

*and*

$$I\left(\frac{1}{2}\right) = \frac{1}{2}(\ln 2)^{\xi+1}, \ \xi \in \mathbb{R}.$$

*(iii) The Shannon entropy of* $X_{(n)}$ *and the similar one of (17) are, respectively,*

$$\mathcal{H}(X_{(n)}) = 1 + \gamma + \xi\gamma + \xi \ln n, \ \xi \in \mathbb{R},$$
$$\mathcal{H}(G_\xi) = 1 + \gamma + \xi\gamma, \ \xi \in \mathbb{R},$$

*while the upper bound* $\mathcal{UB}_{EVD}^{\mathcal{H}}$ *of the extreme value distribution (15) is given by,*

$$\mathcal{UB}_{EVD}^{\mathcal{H}} = 1 - \ln\left[2I\left(\frac{1}{2}\right)\right] + \gamma = 1 + \gamma - (\xi + 1)\ln(\ln 2),$$

*for* $\xi = 0$ *or* $-1 < \xi < 0$, *where* $\gamma$ *is Euler's constant.*
*(iv) The extropy of* $X_{(n)}$ *and the similar one of (17) are, respectively,*

$$\mathcal{J}(X_{(n)}) = -\frac{\Gamma(\xi + 2)}{2^{\xi+3}n^\xi}, \ \xi > -2,$$
$$\mathcal{J}(G_\xi) = -\frac{\Gamma(\xi + 2)}{2^{\xi+3}}, \ \xi > -2,$$

*while the upper bound $\mathcal{UB}_{EVD}^{\mathcal{J}}$ of the extreme value distribution (15) is given by,*

$$\mathcal{UB}_{EVD}^{\mathcal{J}} = -\frac{1}{4} I \left( \frac{1}{2} \right) = -\frac{1}{8} (\ln 2)^{\xi+1},$$

*for $\xi = 0$ or $-1 < \xi < 0$, where $\Gamma$ is the gamma function.*

**Remark 3** **(a)** *Part (i) of the proposition investigates the log-concavity of the generalized extreme value distribution $g_\xi$ in (17). The log-concavity of the case $\xi = 0$ is also included in Table 1 in Bagnoli and Bergstrom (2005). The density $g_\xi$ in (17) is log-concave for a specific range of values of the shape parameter $\xi$, namely, only in the case $\xi = 0$ or when $-1 < \xi < 0$. This should not surprise by taking into account the discussion about a number of non-regular situations associated with $\xi$, which is provided in p. 62 of Kotz and Nadarajah (2000). Based on this discussion, "the experience with real-world data suggests that the condition $-1/2 < \xi < 1/2$ is almost always satisfied in practical applications - in particular in environmetrics." This sentence is also mentioned on p. 1442 of Müller and Rufibach (2008) to underline that the restriction of $\xi = 0$ or $-1 < \xi < 0$ is less restrictive than it seems at first sight. Hence, log-concavity of $g_\xi$ is obeyed for values of $\xi$ which are in harmony with applications of this model, in practice.*
**(b)** *The explicit expression of $\mathcal{H}(G_\xi)$ is also provided by formula (19) in Zaid et al. (2022) and in Saeb (2023). The formula for $\mathcal{J}(G_\xi)$ does not appear in the bibliography to the best of our knowledge.*
**(c)** *Let's now consider the convergence which is formulated in (14) and let's concentrate on the norming constants $a_n$ and $b_n$, $n \geq 1$, with $a_n > 0$. Based on (14) and in view of Corollary 1.2.4 and Remark 1.2.7, p. 9-10, 21-22 of de Haan and Ferreira (2006), or Proposition 1.1, p. 40, of Resnick (2013), or Kotz and Nadarajah (2000), p. 8-9 and Beirland et al. (2004), p. 48, (2.2), the extreme value distribution $G_\xi$ of (15), must be of one of the following three types:*
$G_{1,\xi}(x) = \exp\left(-x^{-1/\xi}\right), x > 0,\ \xi > 0$ *(Fréchet) with $a_n = U(n)$, with $U(t) = F^{-1}\left(1 - \frac{1}{t}\right)$ and $b_n = 0$,*
$G_{2,\xi}(x) = \exp\left(-(-x)^{-1/\xi}\right), x \leq 0, \xi < 0$ *(reversed Weibull) with $a_n = x^* - U(n) = F^{-1}(1) - F^{-1}\left(1 - \frac{1}{n}\right)$ and $b_n = x^* = F^{-1}(1), x^* = \sup\{x : F(x) < 1\}$,*
$G_{3,\xi}(x) = \exp\left(-e^{-x}\right), x \in \mathbb{R}, \xi = 0$ *(Gumbel distribution) with $a_n = h(U(n)) = h\left(F^{-1}\left(1 - \frac{1}{n}\right)\right)$, with $h(u) = \frac{1-F}{F'}$ and $b_n = F^{-1}(1) - F^{-1}\left(1 - \frac{1}{n}\right)$,*
*where $F^{-1}(1) = x^* = \sup\{x : F(x) < 1\}$, taking into account Bobkov and Ledoux (2019), p. 83. The limiting behavior of the function $U(n) = F^{-1}\left(1 - \frac{1}{n}\right)$ is investigated in Lemma 1.2.9 in p. 22 of de Haan and Ferreira (2006) and $\lim_{t \to \infty} U(t)$ is finite, $\lim_{t \to \infty} U(t) < \infty$, if $\xi < 0$.*

The next example, motivated by the Example 1.2 in Johnson (2024), helps to illustrate the evaluation of the norming constant $a_n$, $n \geq 1$, with $a_n > 0$, in the case of exponential distribution.

**Example 1** *Assume that a random sample $X_1, X_2, ..., X_n$ of size $n$ is coming from an exponential distribution function $F_\theta(x) = 1 - e^{-\theta x}$, $x > 0$. Then, from Table 2.3, p. 72, of Beirland et al. (2004), the exponential distribution is in the Gumbel domain with distribution function $G_0(x) = \exp\left(-e^{-x}\right), x \in \mathbb{R}, \xi = 0$, according to (16). In this frame, $F_\theta^{-1}(y) = -\frac{1}{\theta} \ln(1 - y)$ and elementary algebraic manipulations lead to*

$$a_n = h\left(F_\theta^{-1}\left(1 - \frac{1}{n}\right)\right) = \frac{1 - F_\theta\left(F_\theta^{-1}\left(1 - \frac{1}{n}\right)\right)}{\theta \exp\left(-\theta F_\theta^{-1}\left(1 - \frac{1}{n}\right)\right)} = \frac{1}{\theta}.$$

*Conversely, if $a_n = \frac{1}{\theta}$, then, $\frac{1 - F_\theta\left(F_\theta^{-1}\left(1 - \frac{1}{n}\right)\right)}{\theta \exp\left(-\theta F_\theta^{-1}\left(1 - \frac{1}{n}\right)\right)} = \frac{1}{\theta}$ or $\exp\left(-\theta F_\theta^{-1}(1 - y)\right) = y$, for $y = \frac{1}{n} \in (0, 1)$.*
*The last equality entails that $F_\theta^{-1}(z) = -\frac{1}{\theta}\ln(1 - z)$, for $z = 1 - y > 0$ and then, $F_\theta(w) = 1 - e^{-\theta w}$, $w > 0$, that is exponential distribution with parameter $\theta > 0$.*

Having closing the above short exposition of extreme value distributions, for reasons posed at p. 46-47 of Beirlant et al. (2004), lets concentrate on the limiting behavior of Shannon entropy $\mathcal{H}\left(\frac{X_{(n)} - b_n}{a_n}\right)$ and extropy $\mathcal{J}\left(\frac{X_{(n)} - b_n}{a_n}\right)$ as $n \to \infty$, for sequences $\{a_n; n \geq 1\}$ and $\{b_n; n \geq 1\}$ with the first a sequence of positive numbers. Based on Beirlant et al. (2004), p. 46, the standardization with $b_n$ and $a_n$ is natural since otherwise $X_{(n)} \to x_*$. On the other hand, it is well known from information theoretic derivations that the Shannon entropy and the extropy are location-free but scale-dependent, in the sense

$$\mathcal{H}(cX + d) = \ln c + \mathcal{H}(X),$$
$$\mathcal{J}(cX + d) = \frac{1}{c}\mathcal{J}(X), \quad \text{with } c, d \in \mathbb{R}, c > 0,$$

in view of Theorems 8.6.3 and 8.6.4 in Cover and Thomas (2006), p. 253-254 and Tommaj et al. (2023), p. 1337.

Hence, taking into account the above mentioned location-free but scale-dependent behavior of $\mathcal{H}$ and $\mathcal{J}$, it is immediate to see that

$$\mathcal{H}\left(\frac{X_{(n)} - b_n}{a_n}\right) = -\ln a_n + \mathcal{H}(X_{(n)}), \tag{19}$$

and therefore $\lim_{n \to \infty} \mathcal{H}(X_{(n)})$ is directly connected with the limit, as $n \to \infty$, of the left hand side of (19) and the similar limit of the sequence $a_n$, as well. On the other hand, (5) and (19) lead to

$$\mathcal{H}\left(\frac{X_{(n)} - b_n}{a_n}\right) = 1 - \ln a_n - \ln n - \frac{1}{n} - E_Y[\ln I(Y)), \tag{20}$$

where $I(Y)$, $Y \sim Beta(n, 1)$ is defined by (4). In a similar manner, the extropy $\mathcal{J}$ satisfies

$$\mathcal{J}\left(\frac{X_{(n)} - b_n}{a_n}\right) = a_n \mathcal{J}(X_{(n)}), \tag{21}$$

for sequences $\{a_n; n \geq 1\}$ and $\{b_n; n \geq 1\}$ with the first a sequence of positive numbers and $\mathcal{J}(X_{(n)})$ given by (12).

In the setting presented above, the next proposition states bounds for the Shannon entropy and the extropy of the normalizing sequence $(X_{(n)} - b_n)/a_n$.

**Proposition 2** *Let a random sample $X_1, X_2, ..., X_n$ of size $n$ from a continuous probability distribution function $F$ with respective probability density function $f$, supported on some interval $(\alpha, \beta)$, finite or not, where $f$ is positive and $\ln f$ is concave. Then,*
*(i) The following inequalities are true,*

$$1 - \ln n - \frac{1}{n} - \ln a_n \leq \mathcal{H}\left(\frac{X_{(n)} - b_n}{a_n}\right) \leq 1 - \ln\left[2I\left(\frac{1}{2}\right)\right] - \ln n - \frac{1}{n} + \ln 2 - \ln a_n + \sum_{k=1}^{n}\frac{1}{k} - \sum_{k=1}^{n-1}\frac{1}{k 2^k},$$

*and*

$$1 - \lim_{n \to \infty} \left(\ln(na_n)\right) \le \lim_{n \to \infty} \mathcal{H}\left(\frac{X_{(n)} - b_n}{a_n}\right) \le 1 - \ln\left[2I\left(\frac{1}{2}\right)\right] + \gamma - \lim_{n \to \infty} \left(\ln a_n\right).$$

*(ii) It is true that,*

$$-\frac{na_n}{2} I\left(\frac{1}{2}\right) \le \mathcal{J}\left(\frac{X_{(n)} - b_n}{a_n}\right) = a_n \mathcal{J}(X_{(n)}) \le -n^2 a_n I\left(\frac{1}{2}\right)\left(\frac{1}{2n-1} - \frac{1}{2n} - \frac{1}{(2n-1)2^{2n-1}} + \frac{2}{2n2^{2n}}\right),$$

*and*

$$-\frac{1}{2} I\left(\frac{1}{2}\right) \lim_{n \to \infty} (na_n) \le \lim_{n \to \infty} \mathcal{J}\left(\frac{X_{(n)} - b_n}{a_n}\right) \le -\frac{1}{4} I\left(\frac{1}{2}\right) \lim_{n \to \infty} (a_n).$$

The proof of part (i) is immediately obtained from parts (i) and (ii) of Theorem 1 by taking into account (19). The proof of part (ii) is obtained from parts (i) and (ii) of Theorem 3 by taking into account (21).

Motivated from Theorem 2, let's now try to investigate attainment of the upper bound for the $\lim_{n \to \infty} \mathcal{H}\left(\frac{X_{(n)} - b_n}{a_n}\right)$, formulated in part (i) of the above proposition, on the basis of a random sample $X_1, X_2, ..., X_n$ of size $n$ from a continuous probability distribution function $F_\theta$ with respective probability density function $f_\theta$, supported on some interval $(\alpha, \beta)$, finite or not, where $f_\theta$ is positive and $\ln f_\theta$ is concave and they depend on a real parameter $\theta$. In this context, based on (19), (10) and part (i) of Proposition 2,

$$\lim_{n \to \infty} \mathcal{H}\left(\frac{X_{(n)} - b_n}{a_n}\right) \le \lim_{n \to \infty} (-\ln a_n) + \mathcal{UB}^{\mathcal{H}}, \tag{22}$$

with equality if and only if $F_\theta$ is the distribution function of the exponential distribution, in view of Theorem 2. Based on (22), on the previous example where $a_n = \frac{1}{\theta}$ for the exponential distribution and on the fact that $\mathcal{UB}^{\mathcal{H}}_{\exp} = 1 - \ln \theta + \gamma$, we conclude that, equality is valid in (22) if and only if

$$\begin{aligned}\lim_{n \to \infty} \mathcal{H}\left(\frac{X_{(n)} - b_n}{a_n}\right) &= \ln \theta + 1 - \ln \theta + \gamma \\ &= 1 + \gamma.\end{aligned}$$

However, part (iii) of Proposition 1 entails that $1 + \gamma$ is the Shannon entropy $\mathcal{H}(G_0)$ of the Gumbel distribution with distribution function $G_0(x) = \exp\left(-e^{-x}\right), x \in \mathbb{R}$ for $\xi = 0$.

A quite similar analysis for extropy entails that,

$$\lim_{n \to \infty} \mathcal{J}\left(\frac{X_{(n)} - b_n}{a_n}\right) \le \lim_{n \to \infty} (a_n) \mathcal{UB}^{\mathcal{J}}, \tag{23}$$

with equality if and only if the underlined distribution is the exponential distribution. Hence, equality is valid in (23) if and only if

$$\begin{aligned}\lim_{n \to \infty} \mathcal{J}\left(\frac{X_{(n)} - b_n}{a_n}\right) &= \lim_{n \to \infty} \left(\frac{1}{\theta}\right) \times \left(-\frac{\theta}{8}\right) \\ &= -\frac{\theta}{8} = \mathcal{J}(G_0).\end{aligned}$$

Based on the above analysis we can now state the next characterization of attainment of the upper bounds of Proposition 2.

**Theorem 5** *Let a random sample $X_1, X_2, ..., X_n$ of size $n$ from a continuous probability distribution function $F_\theta$ with respective probability density function $f_\theta$, supported on some interval $(\alpha, \beta)$, finite or not, where $f_\theta$ is positive and $\ln f_\theta$ is concave and they depend on a real parameter $\theta$. Then,*

*(i) $\lim_{n\to\infty} \mathcal{H}\left(\frac{X_{(n)} - b_n}{a_n}\right)$ attains its maximum value $\lim_{n\to\infty}(-\ln a_n) + \mathcal{UB}^{\mathcal{H}}$, for $\mathcal{UB}^{\mathcal{H}}$ given in (10), if and only if only $F_\theta$ is the distribution function of the exponential distribution. Moreover,*

$$\lim_{n\to\infty} \mathcal{H}\left(\frac{X_{(n)} - b_n}{a_n}\right) = 1 + \gamma = \mathcal{H}(G_0). \tag{24}$$

*(ii) $\lim_{n\to\infty} \mathcal{J}\left(\frac{X_{(n)} - b_n}{a_n}\right)$ attains its maximum value $\lim_{n\to\infty}(a_n)\mathcal{UB}^{\mathcal{J}}$, for $\mathcal{UB}^{\mathcal{J}}$ given in (13), if and only if only $F_\theta$ is the distribution function of the exponential distribution. Moreover,*

$$\lim_{n\to\infty} \mathcal{J}\left(\frac{X_{(n)} - b_n}{a_n}\right) = -\frac{\theta}{8} = \mathcal{J}(G_0). \tag{25}$$

The above theorem implicitly derives the convergence, presented in Johnson (2024), Ravi and Saeb (2014) and Saeb (2023), of Shannon entropy $\mathcal{H}\left((X_{(n)} - b_n)/a_n\right)$ to the Shannon entropy of the Gumber distribution $G_0$. Moreover, it implicitly extends this convergence for the extropy of the normalizing sequence $(X_{(n)} - b_n)/a_n$, formulated in (25), which has not been appeared in the respective literature to the best of our knowledge. The main assumption in Theorem 5 is that of the log-concavity of the underlined distribution that governs the data. The respective assumption in Ravi and Saeb (2014) and Saeb (2023) to prove the convergence in (24) is that of monotonicity (non-increasing) of the distribution which is considered to govern the data. Johnson (2024) proves the convergence in (24) subject to the assumption that the respective distribution has a von Mises representation, presented in (4.1) of Johnson (2024).

# 4    Conclusions

In this paper, two information theoretic measures are considered, the Shannon entropy and the extropy, and they are applied to the distribution of the largest order statistics and its normalized version in the extreme value theory setting. Bounds of these measures are derived and a characterization of the exponential distribution is provided as the model that maximizes the uncertainty which is associated with the maximum value in a large size random sample from a log-concave distribution. It was moreover proved that the Shannon entropy and the extropy of the normalized maxima of a large size sample convergence to the respective measures for the Gumbel distribution, something which was occupied the recent literature (cf. Saeb (2023), Johnson (2024)). The results derived in this paper have tried to give more insight in the behavior of two omnipresent information theoretic quantities in the extreme value theory regime.

# 5    Appendix: Proofs of the theoretical results

## 5.1    Proof of Theorem 1

(i) Based on the fact that $f$ is log-concave, $\ln f$ is concave and then, on the basis of Bobkov and Madiman (2011), p. 1532, the function $I(t) = f(F^{-1}(t))$, $0 < t < 1$, defined by (7), is positive and concave on $(0,1)$ (cf. also Bobkov (1996), Proposition A.1). Given that $I(t)$ is positive and concave, $\ln I(t)$ is also concave and then Jensen's inequality leads to,

$$E_Y[\ln I(Y)] \le \ln I\left(E_Y(Y)\right), \text{ with } Y \sim Beta(n,1),$$

and therefore

$$E_Y[\ln I(Y)] \le \ln I\left(\frac{n}{n+1}\right), \tag{26}$$

because $E_Y(Y) = n/(n+1)$, for $Y \sim Beta(n,1)$.

On the other hand, taking into account Bobkov and Madiman (2011), p. 1532,

$$\min\{t, 1-t\} \le I(t) \le 1, \; 0 < t < 1,$$

which leads to

$$-\ln(n+1) \le \ln I\left(\frac{n}{n+1}\right) \le 0. \tag{27}$$

Inequalities (26) and (27) entail that

$$-E_Y[\ln I(Y))] \ge -\ln I\left(\frac{n}{n+1}\right) \ge 0,$$

which in association with (5) gives that,

$$\mathcal{H}(X_{(n)}) \ge 1 - \ln n - \frac{1}{n}. \tag{28}$$

Based again on Bobkov (1999), p. 1916, inequality (4.7),

$$2I\left(\frac{1}{2}\right)\min\{t, 1-t\} \le I(t), \; 0 < t < 1,$$

and then,

$$2I\left(\frac{1}{2}\right)\min\{Y, 1-Y\} \le I(Y), \; Y \sim Beta(n,1),$$

or

$$\ln\left[2I\left(\frac{1}{2}\right)\right] + E_Y[\ln(\min\{Y, 1-Y\})] \le E_Y[\ln I(Y)], \; Y \sim Beta(n,1). \tag{29}$$

However,

$$
\begin{aligned}
E_Y[\ln(\min\{Y, 1-Y\})] &= \int_0^1 \ln(\min\{y, 1-y\}) f_{Beta(n,1)}(y) dy = \int_0^1 \ln(\min\{y, 1-y\}) n y^{n-1} dy \\
&= \int_0^{1/2} \ln(\min\{y, 1-y\}) n y^{n-1} dy + \int_{1/2}^1 \ln(\min\{y, 1-y\}) n y^{n-1} dy \\
&= \int_0^{1/2} \ln(y) n y^{n-1} dy + \int_{1/2}^1 \ln(1-y) n y^{n-1} dy.
\end{aligned}
\tag{30}
$$

Based on Gradshteyn and Ryzhik (2007), p. 238, 2.723(1),

$$\int\limits_0^{1/2} \ln(y) n y^{n-1} dy = -\frac{1}{2^n}\left(\ln 2 + \frac{1}{n}\right). \tag{31}$$

Let's now work the second integral. Based on Gradshteyn and Ryzhik (2007), p. 239, 2.728(1),

$$\int\limits_{1/2}^{1} \ln(1-y) n y^{n-1} dy = \left[y^n \ln(1-y) + \int \frac{y^n}{1-y} dy\right]_{1/2}^{1}. \tag{32}$$

Again, using 2.111, $3^8$ of p. 69 in Gradshteyn and Ryzhik (2007),

$$
\begin{aligned}
\int \frac{y^n}{1-y} dy &= -\frac{y^n}{n} - \frac{y^{n-1}}{n-1} - \frac{y^{n-2}}{n-2} - ... + (-1)^{n-1}\frac{y}{1\cdot(-1)^n} + \frac{(-1)^n}{(-1)^{n+1}}\ln(1-y) \\
&= -\left(\frac{y^n}{n} + \frac{y^{n-1}}{n-1} + \frac{y^{n-2}}{n-2} + ... + \frac{y}{1} + \ln(1-y)\right).
\end{aligned} \tag{33}
$$

In this frame, (32) and (33) entail that,

$$
\begin{aligned}
\int\limits_{1/2}^{1} \ln(1-y) n y^{n-1} dy &= \left[y^n \ln(1-y) - \frac{y^n}{n} - \frac{y^{n-1}}{n-1} - \frac{y^{n-2}}{n-2} - ... - \frac{y}{1} - \ln(1-y)\right]_{1/2}^{1} \\
&= \left[[y^n \ln(1-y)]|_{y=1} - [\ln(1-y)]|_{y=1} - \left(\frac{1}{n} + \frac{1}{n-1} + \frac{1}{n-2} + ... + \frac{1}{1}\right)\right] \\
&\quad - \left[-\left(\frac{1}{2^n} - 1\right)\ln 2 - \left(\frac{1}{2^n n} + \frac{1}{2^{n-1}(n-1)} + \frac{1}{2^{n-2}(n-2)} + ... + \frac{1}{2\cdot1}\right)\right] \\
&= -\sum_{k=1}^{n} \frac{1}{k} + \sum_{k=1}^{n}\frac{1}{k2^k} + \left(\frac{1}{2^n} - 1\right)\ln 2.
\end{aligned} \tag{34}
$$

Based now in (30), (31) and (34),

$$E_Y[\ln(\min\{Y, 1-Y\})] = -\frac{1}{2^n}\left(\ln 2 + \frac{1}{n}\right) - \sum_{k=1}^{n}\frac{1}{k} + \sum_{k=1}^{n}\frac{1}{k2^k} + \left(\frac{1}{2^n} - 1\right)\ln 2,$$

or

$$
\begin{aligned}
E_Y[\ln(\min\{Y, 1-Y\})] &= -\frac{1}{2^n}\ln 2 - \frac{1}{n2^n} - \sum_{k=1}^{n}\frac{1}{k} + \sum_{k=1}^{n}\frac{1}{k2^k} + \frac{1}{2^n}\ln 2 - \ln 2 \\
&= -\ln 2 - \sum_{k=1}^{n}\frac{1}{k} + \sum_{k=1}^{n-1}\frac{1}{k2^k}.
\end{aligned} \tag{35}
$$

Then, based on (29) and (35),

$$\ln\left[2I\left(\frac{1}{2}\right)\right] + E_Y[\ln(\min\{Y, 1-Y\})] = \ln\left[2I\left(\frac{1}{2}\right)\right] - \ln 2 - \sum_{k=1}^{n}\frac{1}{k} + \sum_{k=1}^{n-1}\frac{1}{k2^k} \le E_Y[\ln I(Y)],$$

with $Y \sim Beta(n, 1)$ and

$$-E_Y[\ln I(Y)] \leq -\ln\left[2I\left(\frac{1}{2}\right)\right] + \ln 2 + \sum_{k=1}^{n}\frac{1}{k} - \sum_{k=1}^{n-1}\frac{1}{k2^k}, \ Y \sim Beta(n, 1). \tag{36}$$

Hence, (5) and (36) give that

$$\begin{aligned}
\mathcal{H}(X_{(n)}) &= 1 - \ln n - \frac{1}{n} - E_Y[\ln I(Y)) \\
&\leq 1 - \ln\left[2I\left(\frac{1}{2}\right)\right] - \ln n - \frac{1}{n} + \ln 2 + \sum_{k=1}^{n}\frac{1}{k} - \sum_{k=1}^{n-1}\frac{1}{k2^k}. \tag{37}
\end{aligned}$$

In a summary, if $f$ is log-concave, then from (28) and (37) we conclude that,

$$1 - \ln n - \frac{1}{n} \leq \mathcal{H}(X_{(n)}) \leq 1 - \ln\left[2I\left(\frac{1}{2}\right)\right] - \ln n - \frac{1}{n} + \ln 2 + \sum_{k=1}^{n}\frac{1}{k} - \sum_{k=1}^{n-1}\frac{1}{k2^k}. \tag{38}$$

(ii) It is well known (cf., among others, Abramowitz and Stegun (1970), p. 255, 258, Mathai (1993), p. 4, 11) that

$$\gamma = \lim_{n\to\infty}(H_n - \ln n) = \lim_{n\to\infty}\left(\sum_{k=1}^{n}\frac{1}{k} - \ln n\right), \tag{39}$$

where $\gamma$ is the Euler constant while $H_n$ is the $n$-th harmonic number, defined by $H_n = \sum_{k=1}^{n}\frac{1}{k}$.
Then, based on part (i) of the proposition, it remains to obtain the $\lim_{n\to\infty}\sum_{k=1}^{n-1}\frac{1}{k2^k}$. Working in this direction and based on the sum of an infinite geometric series, it is true that,

$$\sum_{k=1}^{\infty}\left(\frac{z}{2}\right)^{k-1} = \frac{1}{1-(z/2)} = \frac{2}{2-z}, \ |z| < 2.$$

Then, integrating both sides

$$\sum_{k=1}^{\infty}\frac{1}{2^{k-1}}\int_0^1\left(\frac{z^k}{k}\right)' dz = 2\int_0^1\frac{1}{2-z}dz,$$

or

$$\sum_{k=1}^{\infty}\frac{1}{2^{k-1}}\frac{1}{k} = -2\int_0^1[\ln(2-z)]'dz,$$

which gives that

$$\sum_{k=1}^{\infty}\frac{1}{2^{k-1}}\frac{1}{k} = 2\ln 2$$

and therefore

$$\sum_{k=1}^{\infty}\frac{1}{k2^k} = \ln 2. \tag{40}$$

Then, part (i) of the proposition along with (39) and (40) complete the proof of part (ii). ▲

## 5.2  Proof of Theorem 2

($\Longrightarrow$) Let the random sample is coming from the exponential distribution with density $f_\theta(x) = \theta e^{-\theta x}$, $x > 0$ and parameter $\theta > 0$. Then, from the Table 1, the assertion of the theorem is valid.

($\Longleftarrow$) Let $\lim_{n\to\infty} \mathcal{H}(X_{(n)}) = \mathcal{UB}^\mathcal{H}$. Then, from (6) and (10),

$$\lim_{n\to\infty} \mathcal{H}(X_{(n)}) = \lim_{n\to\infty} \left( 1 - \ln n - \frac{1}{n} - \int_0^1 n y^{n-1} \ln f_\theta(F_\theta^{-1}(y)) dy \right) = 1 - \ln\left[ 2I_\theta\left(\frac{1}{2}\right) \right] + \gamma.$$

Taking into account (8), $\lim_{n\to\infty} \ln n = \lim_{n\to\infty} H_n - \gamma$ and hence $\lim_{n\to\infty} \mathcal{H}(X_{(n)}) = \mathcal{UB}^\mathcal{H}$ entails that,

$$1 - \lim_{n\to\infty} (\ln n) - \lim_{n\to\infty} \frac{1}{n} - \lim_{n\to\infty} \int_0^1 n y^{n-1} \ln I_\theta(y)) dy = 1 - \ln\left[ 2I_\theta\left(\frac{1}{2}\right) \right] + \gamma,$$

or

$$\lim_{n\to\infty} H_n + \lim_{n\to\infty} \int_0^1 n y^{n-1} \ln I_\theta(y) dy = \ln\left[ 2I_\theta\left(\frac{1}{2}\right) \right],$$

or, taking into account that $n y^{n-1}, 0 < y < 1$, is the density of $Beta(n,1)$, it is valid that

$$\int_0^1 n y^{n-1} \ln\left[ 2I_\theta\left(\frac{1}{2}\right) \right] dy = \ln\left[ 2I_\theta\left(\frac{1}{2}\right) \right],$$

and therefore the equation $\lim_{n\to\infty} \mathcal{H}(X_{(n)}) = \mathcal{UB}^\mathcal{H}$ entails that,

$$\lim_{n\to\infty} \left\{ H_n + \int_0^1 n y^{n-1} \ln \frac{I_\theta(y)}{2I_\theta(1/2)} dy \right\} = 0, \tag{41}$$

where,

$$I_\theta(y) = f_\theta(F_\theta^{-1}(y)), 0 < y < 1,$$

in view of (7). Moreover, in view of well known properties of the beta distribution (cf. Zografos and Balakrishnan (2009), Lemma 1(c) and Nawa and Nadarajah (2023), p.2, equation (3)), for $Y \sim Beta(n,1)$,

$$E[\ln(1 - Y)] = \Psi(1) - \Psi(n+1) = -H_n,$$

and therefore

$$H_n = -E[\ln(1 - Y)] = -\int_0^1 n y^{n-1} \ln(1 - y) dy. \tag{42}$$

Based on (41) and (42), $\lim_{n\to\infty} \mathcal{H}(X_{(n)}) = \mathcal{UB}^\mathcal{H}$ entails that,

$$\lim_{n\to\infty} \int_0^1 n y^{n-1} \ln \frac{I_\theta(y)}{2I_\theta(1/2)(1-y)} dy = 0. \tag{43}$$

The last equation can be written as follows,

$$\lim_{n\to\infty}\int\limits_0^1 ny^{n-1}\ln\frac{I_\theta(y)}{2I_\theta(1/2)\min\{y,1-y\}}dy + \lim_{n\to\infty}\int\limits_0^1 ny^{n-1}\ln\frac{\min\{y,1-y\}}{1-y}dy = 0. \qquad (44)$$

Consider now the right hand side integral of (44) which is simplified as follows,

$$
\begin{aligned}
\int\limits_0^1 ny^{n-1}\ln\frac{\min\{y,1-y\}}{1-y}dy &= \int\limits_0^{1/2} ny^{n-1}\ln\frac{y}{1-y}dy + \int\limits_{1/2}^1 ny^{n-1}\ln\frac{1-y}{1-y}dy \\
&= \int\limits_0^{1/2} ny^{n-1}\ln y\,dy - \int\limits_0^{1/2} ny^{n-1}\ln(1-y)dy. \qquad (45)
\end{aligned}
$$

Based on 2.723(1) of p. 238 in Gradshteyn and Ryzhik (2007) or on the proof of Theorem 1,

$$\int\limits_0^{1/2} ny^{n-1}\ln y\,dy = -\frac{1}{2^n}\ln 2 - \frac{1}{n2^n}. \qquad (46)$$

Moreover, based on 2.728(1) of p. 239 in Gradshteyn and Ryzhik (2007) or on the proof of Theorem 1,

$$
\begin{aligned}
\int\limits_0^{1/2} ny^{n-1}\ln(1-y)dy &= \left[y^n\ln(1-y) + \int\frac{y^n}{1-y}dy\right]_0^{1/2} \\
&= \frac{1}{2^n}\ln\frac{1}{2} + \left[\int\frac{y^n}{1-y}dy\right]_0^{1/2}. \qquad (47)
\end{aligned}
$$

However, from 2.111, $3^8$ of p. 69 in Gradshteyn and Ryzhik (2007) or (33)

$$\int\frac{y^n}{1-y}dy = -\left(\frac{y^n}{n} + \frac{y^{n-1}}{n-1} + \frac{y^{n-2}}{n-2} + ... + \frac{y}{1} + \ln(1-y)\right)$$

and therefore

$$\left[\int\frac{y^n}{1-y}dy\right]_0^{1/2} = -\left(\frac{1}{n2^n} + \frac{1}{(n-1)2^{n-1}} + ... + \frac{1}{2} + \ln\frac{1}{2}\right). \qquad (48)$$

Based now in (46), (47) and (48),

$$\int\limits_0^1 ny^{n-1}\ln\frac{\min\{y,1-y\}}{1-y}dy = \frac{1}{(n-1)2^{n-1}} + \frac{1}{(n-2)2^{n-2}} + ... + \frac{1}{2} + \ln\frac{1}{2},$$

or

$$\int\limits_0^1 ny^{n-1}\ln\frac{\min\{y,1-y\}}{1-y}dy = \sum_{k=1}^{n-1}\frac{1}{k2^k} - \ln 2. \qquad (49)$$

However, in the proof of Theorem 1 has been shown (cf. equation (40)) that,

$$\ln 2 = \sum_{k=1}^{\infty} \frac{1}{k 2^k}. \tag{50}$$

Then, (49) and (50) lead to

$$\int\limits_0^1 n y^{n-1} \ln \frac{\min\{y, 1-y\}}{1-y} dy = \sum_{k=1}^{n-1} \frac{1}{k 2^k} - \sum_{k=1}^{\infty} \frac{1}{k 2^k}$$

and therefore

$$\lim_{n \to \infty} \int\limits_0^1 n y^{n-1} \ln \frac{\min\{y, 1-y\}}{1-y} dy = \lim_{n \to \infty} \sum_{k=1}^{n-1} \frac{1}{k 2^k} - \sum_{k=1}^{\infty} \frac{1}{k 2^k} = 0. \tag{51}$$

Summarizing all the above, from (43), (44) and (51),

$$\lim_{n \to \infty} \mathcal{H}(X_{(n)}) = \mathcal{UB}^{\mathcal{H}} \text{ entails that } \lim_{n \to \infty} \int\limits_0^1 n y^{n-1} \ln \frac{I_\theta(y)}{2 I_\theta (1/2) \min\{y, 1-y\}} dy = 0. \tag{52}$$

This last formula (52) means that $\lim_{n \to \infty} \mathcal{H}(X_{(n)}) = \mathcal{UB}^{\mathcal{H}}$ entails that

$$\int\limits_0^1 n y^{n-1} \ln \frac{I_\theta(y)}{2 I_\theta (1/2) \min\{y, 1-y\}} dy = 0,$$

for large $n$ ($n \to \infty$). Given that $n y^{n-1}, 0 < y < 1$, is the density of $Beta(n, 1)$ distribution, the function inside the last integral is non-negative because $I_\theta(y) \geq 2 I_\theta (1/2) \min\{y, 1-y\}$, in view of (4.7) of p. 1916 in Bobkov (1999). Therefore,

$$\lim_{n \to \infty} \mathcal{H}(X_{(n)}) = \mathcal{UB}^{\mathcal{H}} \text{ entails that } I_\theta(y) = 2 I_\theta (1/2) \min\{y, 1-y\}, \text{ for large } n \ (n \to \infty). \tag{53}$$

On the other hand, an application of (7) for $t = F(y)$ leads to

$$I(F(y)) = f(F^{-1}(F(y))) = f(y). \tag{54}$$

Without any loss of generality suppose that $\min\{y, 1-y\} = 1-y, 0 < y < 1$. Then, (53) and (54) give that,

$$\lim_{n \to \infty} \mathcal{H}(X_{(n)}) = \mathcal{UB}^{\mathcal{H}} \text{ entails that } I_\theta(F_\theta(y)) = 2 I_\theta (1/2) \min\{F_\theta(y), 1-F_\theta(y)\}, \text{ for large } n \ (n \to \infty),$$

or,

$$\lim_{n \to \infty} \mathcal{H}(X_{(n)}) = \mathcal{UB}^{\mathcal{H}} \text{ entails that } f_\theta(y) = 2 I_\theta (1/2) (1-F_\theta(y)), \ 0 < y < 1, \text{ for large } n \ (n \to \infty). \tag{55}$$

But, $f_\theta(y) = 2 I_\theta (1/2) (1-F_\theta(y))$ in (55) means that $F_\theta(y) = 1 - [2 I_\theta (1/2)]^{-1} f_\theta(y)$ or $(d/dy) F_\theta(y) = -[2 I_\theta (1/2)]^{-1} f'_\theta(y)$ or $f'_\theta(y)/f_\theta(y) = -2 I_\theta (1/2)$. Hence,

$$\frac{d}{dy} \ln f_\theta(y) = -2 I_\theta (1/2) \text{ or } \ln f_\theta(y) = -2 I_\theta (1/2) y + c \text{ or } f_\theta(y) = e^c \exp\{-2 I_\theta (1/2) y\}$$

and given that $\int f_\theta(y)dy = 1$, we obtain that $c = 2I_\theta(1/2)$ and then

$$f_\theta(y) = 2I_\theta(1/2)\exp\{-2I_\theta(1/2)y\}.$$

However, if $g(y) = \lambda e^{-\lambda y}$, then $I_\lambda(t) = g(G^{-1}(t)) = \lambda(1-t), 0 < t < 1$ and therefore $2I_\lambda(1/2) = \lambda$. Therefore, taking into account that the exponential distribution is identifiable we finally obtain that $\lim_{n\to\infty}\mathcal{H}(X_{(n)}) = \mathcal{UB}^{\mathcal{H}}$ entails that $F_\theta$, is the distribution function of the exponential distribution, for large $n$ $(n \to \infty)$, and the theorem has been proved. ▲

## 5.3   Proof of Theorem 3

(i) Based on Bobkov (1999), p. 1916,

$$2I\left(\frac{1}{2}\right)\min\{t, 1-t\} \le I(t) \le I\left(\frac{1}{2}\right) + u\left(t - \frac{1}{2}\right),\ 0 < t < 1,\ |u| \le 2I\left(\frac{1}{2}\right).$$

Hence, taking into account Lemma 1,

$$\mathcal{J}(X_{(n)}) = -\frac{n^2}{2}\int_0^1 t^{2n-2}I(t)dt,$$

and the previous inequality is translated to,

$$-\frac{n^2}{2}I\left(\frac{1}{2}\right)\int_0^1 t^{2n-2}dt - \frac{n^2}{2}u\int_0^1 t^{2n-2}\left(t - \frac{1}{2}\right)dt \le \mathcal{J}(X_{(n)}) \le -n^2 I\left(\frac{1}{2}\right)\int_0^1 t^{2n-2}\min\{t, 1-t\}dt. \tag{56}$$

Calculation of the integrals gives,

$$\int_0^1 t^{2n-2}dt = \left.\frac{t^{2n-1}}{2n-1}\right|_0^1 = \frac{1}{2n-1}, \tag{57}$$

$$\int_0^1 t^{2n-2}\left(t - \frac{1}{2}\right)dt = \left.\frac{t^{2n}}{2n}\right|_0^1 - \frac{1}{2}\left.\frac{t^{2n-1}}{2n-1}\right|_0^1 = \frac{1}{2n} - \frac{1}{2}\frac{1}{2n-1} = \frac{1}{2}\frac{n-1}{n(2n-1)}. \tag{58}$$

On the other hand,

$$\begin{aligned}
\int_0^1 t^{2n-2}\min\{t, 1-t\}dt &= \int_0^{1/2} t^{2n-2}tdt + \int_{1/2}^1 t^{2n-2}(1-t)dt \\
&= \left.\frac{t^{2n}}{2n}\right|_0^{1/2} + \left.\frac{t^{2n-1}}{2n-1}\right|_{1/2}^1 - \left.\frac{t^{2n}}{2n}\right|_{1/2}^1 \\
&= \frac{1}{2n2^{2n}} + \frac{1}{2n-1} - \frac{1}{(2n-1)2^{2n-1}} - \frac{1}{2n} + \frac{1}{2n2^{2n}} \\
&= \frac{1}{2n-1} - \frac{1}{(2n-1)2^{2n-1}} - \frac{1}{2n} + \frac{2}{2n2^{2n}}. 
\end{aligned} \tag{59}$$

Hence, based on (56), (57), (58) and (59), we obtain,

$$-\frac{n^2}{2(2n-1)}I\left(\frac{1}{2}\right)-\frac{n(n-1)}{4(2n-1)}u \leq \mathcal{J}(X_{(n)}) \leq -n^2 I\left(\frac{1}{2}\right)\left(\frac{1}{2n-1}-\frac{1}{2n}-\frac{1}{(2n-1)2^{2n-1}}+\frac{2}{2n2^{2n}}\right),$$

and taking into account that $|u| \leq 2I\left(\frac{1}{2}\right)$ (cf. Bobkov (1999), p. 1916), we finally get,

$$-\frac{n^2}{2(2n-1)}I\left(\frac{1}{2}\right)-\frac{n(n-1)}{4(2n-1)}2I\left(\frac{1}{2}\right) \leq \mathcal{J}(X_{(n)}) \leq -n^2 I\left(\frac{1}{2}\right)\left(\frac{1}{2n-1}-\frac{1}{2n}-\frac{1}{(2n-1)2^{2n-1}}+\frac{2}{2n2^{2n}}\right),$$

or

$$-\frac{2n^2-n}{2(2n-1)}I\left(\frac{1}{2}\right) \leq \mathcal{J}(X_{(n)}) \leq -n^2 I\left(\frac{1}{2}\right)\left(\frac{1}{2n-1}-\frac{1}{2n}-\frac{1}{(2n-1)2^{2n-1}}+\frac{2}{2n2^{2n}}\right),$$

or

$$-\frac{n}{2}I\left(\frac{1}{2}\right) \leq \mathcal{J}(X_{(n)}) \leq -n^2 I\left(\frac{1}{2}\right)\left(\frac{1}{2n-1}-\frac{1}{2n}-\frac{1}{(2n-1)2^{2n-1}}+\frac{2}{2n2^{2n}}\right).$$

(ii) The proof follows immediately in view of (i) and the fact that $\lim_{n\to\infty}\left(\frac{1}{(2n-1)2^{2n-1}}\right) = 0$, $\lim_{n\to\infty}\left(\frac{1}{2n2^{2n}}\right) = 0$ and $\lim_{n\to\infty}\left\{-n^2\left(\frac{1}{2n-1}-\frac{1}{2n}\right)\right\} = -\frac{1}{4}$. ▲

## 5.4 Proof of Theorem 4

($\Longrightarrow$) Let the random sample is coming from the exponential distribution with parameter $\theta > 0$. Then, from the Table 2, the assertion of the theorem is valid.

($\Longleftarrow$) Let $\lim_{n\to\infty}\mathcal{J}(X_{(n)}) = \mathcal{UB}^{\mathcal{J}}$. Then, based on Lemma 1,

$$\mathcal{J}(X_{(n)}) = -\frac{n^2}{2}\int_0^1 t^{2n-2}I_\theta(t)dt$$

and based on (13),

$$\lim_{n\to\infty}\mathcal{J}(X_{(n)}) = \mathcal{UB}^{\mathcal{J}} \text{ leads to } \lim_{n\to\infty}\left\{n^2\int_0^1 t^{2n-2}I_\theta(t)dt\right\} = \frac{1}{2}I_\theta\left(\frac{1}{2}\right). \tag{60}$$

Taking into account that the density of $Y \sim Beta(2n-1, 1)$ is given by $(2n-1)y^{2n-2}, 0 < y < 1$, equation (60) entails that

$$\lim_{n\to\infty}\left\{\frac{n^2}{2n-1}\int_0^1 I_\theta(t)(2n-1)t^{2n-2}dt\right\} = \int_0^1 \frac{1}{2}I_\theta\left(\frac{1}{2}\right)(2n-1)t^{2n-2}dt.$$

Hence, $\lim_{n\to\infty}\mathcal{J}(X_{(n)}) = \mathcal{UB}^{\mathcal{J}}$ leads to

$$\lim_{n\to\infty}\left\{\int_0^1\left[I_\theta(t)-\frac{1}{2}\frac{2n-1}{n^2}I_\theta\left(\frac{1}{2}\right)\right]\frac{n^2}{2n-1}(2n-1)t^{2n-2}dt\right\} = 0. \tag{61}$$

In view of (4.7) of p. 1916 in Bobkov (1999), $I_\theta(t) \geq 2I_\theta(1/2)\min\{t, 1-t\}$, $0 < t < 1$. On the other hand, for large $n$, $n \to \infty$,

$$2\min\{t, 1-t\} \geq \frac{1}{2}\frac{2n-1}{n^2}, 0 < t < 1,$$

because the right hand side convergences to zero. Hence, for large $n$, $n \to \infty$,

$$2\min\{t, 1-t\}I_\theta\left(\frac{1}{2}\right) \geq \frac{1}{2}\frac{2n-1}{n^2}I_\theta\left(\frac{1}{2}\right), 0 < t < 1.$$

Therefore, for large $n$, $n \to \infty$, and $0 < t < 1$,

$$I_\theta(t) - \frac{1}{2}\frac{2n-1}{n^2}I_\theta\left(\frac{1}{2}\right) \geq I_\theta(t) - 2\min\{t, 1-t\}I_\theta\left(\frac{1}{2}\right) \geq 0. \tag{62}$$

Taking into account, in view (62), that

$$I_\theta(t) - \frac{1}{2}\frac{2n-1}{n^2}I_\theta\left(\frac{1}{2}\right) \geq 0,$$

(61) means that

$$I_\theta(t) - \frac{1}{2}\frac{2n-1}{n^2}I_\theta\left(\frac{1}{2}\right) = 0, \text{ for large } n, n \to \infty. \tag{63}$$

Formulas (62) and (63), are summarized to

$$I_\theta(t) - 2\min\{t, 1-t\}I_\theta\left(\frac{1}{2}\right) = 0, \text{ for large } n, n \to \infty. \tag{64}$$

Then, based on (60) and (64),

$$\lim_{n\to\infty}\mathcal{J}(X_{(n)}) = \mathcal{UB}^{\mathcal{J}} \text{ entails that } I_\theta(y) = 2I_\theta(1/2)\min\{t, 1-t\}, \text{ for large } n, n \to \infty,$$

and the proof is completed by following the steps of the proof of Theorem 2, just after (53). ▲

## 5.5 Proof of Proposition 1

(i) It is immediate to see that the second derivative

$$\frac{d^2}{dx^2}\ln g_\xi(x) = (\xi+1)(1+\xi x)^{-2}\left[\xi - (1+\xi x)^{-1/\xi}\right],$$

is negative for $-1 < \xi < 0$ and hence, $g_\xi$ in (17) is log-concave for $-1 < \xi < 0$. The case $\xi = 0$ is known from Table 1 in Bagnoli and Bergstrom (2005).

(ii) It is simple to obtain the inverse of the distribution function $G_\xi$ in (15) and it is given by

$$G_\xi^{-1}(t) = -\frac{1}{\xi} + \frac{1}{\xi}\left(-\frac{1}{\ln t}\right)^\xi, \ 0 < t < 1 \text{ and } \xi \neq 0.$$

Then, straightforward manipulations lead to the expression of $I(t) = g_\xi(G_\xi^{-1}(t))$, $0 < t < 1$, by taking into account (17) and $G_\xi^{-1}$.

In a similar manner, if $\xi = 0$ then in view of (16),

$$G_0^{-1}(t) = -\ln(-\ln t), 0 < t < 1,$$

and taking into account (18), we obtain

$$
\begin{aligned}
I(t) &= g_0(G_0^{-1}(t)) = \exp\left(-\left(-\ln(-\ln t) + e^{\ln(-\ln t)}\right)\right), \\
&= -t\ln t, \ 0 < t < 1.
\end{aligned}
$$

(iii) Based on (6), (7) and part (ii) of the proposition,

$$\mathcal{H}(X_{(n)}) = 1 - \ln n - \frac{1}{n} - \int_0^1 ny^{n-1}\ln\left(y(-\ln y)^{\xi+1}\right)dy. \tag{65}$$

However,

$$\int_0^1 ny^{n-1}\ln\left(y(-\ln y)^{\xi+1}\right)dy = n\int_0^1 y^{n-1}\ln y\,dy + n(\xi+1)\int_0^1 y^{n-1}\ln(-\ln y)dy. \tag{66}$$

Based on Gradshteyn and Ryzhik (2007), 4.253(1), p. 540,

$$\int_0^1 x^{\mu-1}\ln x\,dx = B(\mu,1)[\Psi(\mu) - \Psi(\mu+1)], \ \mu > 0,$$

hence,

$$\int_0^1 y^{n-1}\ln y\,dy = B(n,1)[\Psi(n) - \Psi(n+1)] = \frac{\Gamma(n)\Gamma(1)}{\Gamma(n+1)}\left(-\frac{1}{n}\right) = -\frac{1}{n^2}. \tag{67}$$

It remains to obtain the last integral of the right hand side of (66). Based again on 4.325(8), on p. 570 of Gradshteyn and Ryzhik (2007),

$$\int_0^1 x^{\mu-1}\ln\ln\frac{1}{x}dx = -\frac{1}{\mu}(\gamma + \ln\mu), \ \mu > 0$$

and then,

$$\int_0^1 y^{n-1}\ln(-\ln y)dy = -\frac{1}{n}(\gamma + \ln n). \tag{68}$$

From (66), (67) and (68),

$$\int_0^1 ny^{n-1}\ln\left(y(-\ln y)^{\xi+1}\right)dy = -\frac{1}{n} - (\xi+1)(\gamma + \ln n). \tag{69}$$

Then, based on (65) and (69),

$$\mathcal{H}(X_{(n)}) = 1 - \ln n - \frac{1}{n} + \frac{1}{n} + (\xi+1)(\gamma + \ln n),$$

which leads to the desired result. $\mathcal{H}(G_\xi)$ is immediately obtained from $\mathcal{H}(X_{(n)})$ for $n = 1$ because the density $g_\xi$ of $G_\xi$ is obtained from the density of $X_{(n)}$ for $n = 1$. Last, the upper bound $\mathcal{UB}_{EVD}^{\mathcal{H}} = 1 - \ln\left[2I\left(\frac{1}{2}\right)\right] + \gamma$ is directly obtained in view of part (ii) of the proposition.

(iv) Based on Lemma 1,

$$\mathcal{J}(X_{(n)}) = -\frac{n^2}{2}\int_0^1 t^{2n-2}I(t)dt,$$

with $I(t) = t(-\ln t)^{\xi+1}$, $0 < t < 1$, in view of part (ii). Hence,

$$\mathcal{J}(X_{(n)}) = -\frac{n^2}{2}\int_0^1 t^{2n-1}\left(\ln\frac{1}{t}\right)^{\xi+2-1}dt. \tag{70}$$

From 4.272(6), on p. 551 of Gradshteyn and Ryzhik (2007),

$$\int_0^1 x^{\nu-1}\left(\ln\frac{1}{x}\right)^{\mu-1}dx = \frac{1}{\nu^\mu}\Gamma(\mu), \ \ \nu, \mu > 0. \tag{71}$$

Then, from (70) and (71),

$$\mathcal{J}(X_{(n)}) = -\frac{n^2}{2}\frac{1}{(2n)^{\xi+2}}\Gamma(\xi+2), \ \text{for } \xi > -2,$$

which leads to the desired result. $\mathcal{J}(G_\xi)$ is obtained from $\mathcal{J}(X_{(n)})$ for $n = 1$, while the upper bound $\mathcal{UB}_{EVD}^{\mathcal{J}} = -\frac{1}{4}I\left(\frac{1}{2}\right)$ is directly obtained in view of part (ii) of the proposition. ▲